\documentclass[a4paper]{article}
\usepackage{amsmath,amssymb,theorem}
\newtheorem{theorem}{Theorem}
\newtheorem{proposition}{Proposition}
\newtheorem{lemma}{Lemma}
\newtheorem{corollary}{Corollary}

\newcommand{\Z}{\mathbb{Z}}
\newcommand{\F}{\mathbb{F}}
\renewcommand{\O}{{\mathcal{O}}}

\begin{document}

\title{Fast arithmetic in unramified $p$-adic fields}
\author{Hendrik Hubrechts}
\date{June 30, 2009}
\maketitle

\begin{abstract}
Let $p$ be prime and $\Z_{p^n}$ the degree $n$ unramified extension of the ring of $p$-adic integers $\Z_p$. In this paper we give an overview of some very fast algorithms for common operations in $\Z_{p^n}$ modulo $p^N$. Combining existing methods with recent work of Kedlaya and Umans about modular composition of polynomials, we achieve quasi-linear time algorithms in the parameters $n$ and $N$, and quasi-linear or quasi-quadratic time in $\log p$, for most basic operations on these fields, including Galois conjugation, Teichm\"{u}ller lifting and computing minimal polynomials.
\end{abstract}

\section{Introduction}

In this article we aim at explaining very fast methods for arithmetic in finite precision degree $n$ unramified $p$-adic rings $\Z_{p^n}$. Although no truly new computational ideas are presented, most results are new and follow from combining existing algorithms with recent results of Kedlaya and Umans (in particular Theorem \ref{thmKU} below). We do not intend to give complete algorithms, but rather accurate references to the literature combined with precise asymptotic estimates, so that our results can be used as reference.

A central source for classical fast algorithms is the book \cite{MCA}, and for more specific $p$-adic methods we refer to Chapter 12 of \cite{HEHCC}. The result that allows us to give improvements upon the literature is the following. Note that we will give another version of this theorem, more suitable for our needs, below.

\begin{theorem}\label{thmKU}(Kedlaya -- Umans, \cite[Theorem 7.1 with parameters $m=1$, $N=d$]{KU}) Let $R$ be a finite ring of cardinality $q$ given as $(\Z/r\Z)[Z]/(E(Z))$ for some monic polynomial $E(Z)$. For every constant $\delta>0$ there is an algorithm that does the following. Given polynomials $f(X)$, $g(X)$ and $h(X)$ over $R$ of degree at most $d$, such that $h$ has a unit as leading coefficient and that we have access to $d^{1+\delta}$ distinct elements of $R$ whose differences are units in $R$; then it can compute $f(g(X))\bmod h(X)$ in time $\O\left(d^{1+\delta}\log^{1+o(1)}q\right)$.\end{theorem}

All results below for computing in $\Z_{p^n}$ with precision $p^N$ are quasi-linear except for some extra factor $\log p$ arising from computing a $p$-th power in the finite field $\F_{p^n}$. For example, computing a Teichm\"{u}ller lift requires time $\O\left( (Nn\log^2 p)^{1+\epsilon}\right)$, whereas the most general algorithm in \cite{HEHCC} requires time $\O\left( (Nn^2\log^2p)^{1+\epsilon}\right)$. We note that any improvement in computing $x^p$ in $\F_{p^n}\cong \F_p[x]/\bar\varphi(x)$ over the classical complexity $\O\left( (n\log^2p)^{1+\epsilon}\right)$ of repeated squaring would yield a similar improvement for most of our results. Moreover, it is easy to verify that the memory requirements for all results in this paper are essentially linear.

The main application that we have in mind are $p$-adic point counting algorithms on varieties over finite fields, see e.g.\ Chapter 17 in \cite{HEHCC}, which profit a lot from fast $p$-adic arithmetic. For example, in our papers \cite{Hub1} and \cite{Hub2} we showed how to compute the zeta function of hyperelliptic curves in certain families over a finite field $\F_{p^n}$ in time $\O\left(n^{2.667}\right)$ (for fixed genus and characteristic). This improves to $\O\left(n^{2+\epsilon}\right)$ using the results from this paper.

The structure of the sequel of the paper is quite straightforward: we start in each subsection with a precise formulation of the result, and then give references or a proof. The following results are presented:
elementary arithmetic, Newton iteration, Galois conjugation, equations involving the Frobenius automorphism, Teichm\"{u}ller lift, minimal polynomial, trace, norm and Teichm\"{u}ller modulus.

\section{Fast arithmetic}

We choose for once and for all a prime number $p$, an extension degree $n\geq 1$ and a $p$-adic precision $N\geq 1$. We work in the unramified $p$-adic ring $\Z_{p^n}$ modulo $p^N$, and this field is supposed to be represented as $\Z_p[x]/\varphi(x)$ for some monic inert (i.e.\ irreducible modulo $p$) polynomial $\varphi(x)\in\Z_p[x]$ of degree $n$ and precision $p^N$. From now on the notation $\Z_{p^n}\bmod p^N$ will be used for this setting (including the implicit polynomial $\varphi(x)$).

It is not in the scope of this text to discuss how to find a (large) prime $p$ and some inert polynomial $\varphi(x)$ of given degree $n$. However, we note that for finding $\varphi(x)$ it suffices to compute an irreducible polynomial $\bar\varphi(x)$ of degree $n$ over $\F_p$, which is an extensively studied problem, see e.g.\ the reference in the proof of Theorem \ref{thmKU1} below.

For our purposes Theorem \ref{thmKU} is not immediately applicable, hence we give a reformulation.

\begin{theorem}\label{thmKU1}Let $f(x)$, $g(x)$ and $h(x)$ be polynomials of degree at most $n$ over $\Z_p[x]\bmod p^N$, with $f(x)$ monic. Then we can compute $f(g(x))\bmod h(x)$ in time $\O\left( (Nn\log p)^{1+\epsilon}\right)$.
\end{theorem}
\textsc{Proof.} If $p$ is large enough, say $p\geq n^2$, we can use Theorem \ref{thmKU} directly because $\Z_p$ contains enough (readily available) units. Suppose hence $p<n^2$. Shoup has shown in \cite{ShFastConstr} how to construct some irreducible polynomial $\bar E(Y)$ over $\F_p$ of degree $a$ in time $\O\left(((a^2+a\log p)\log p)^{1+\epsilon}\right)$, subsequently improved in Section 8.5 of \cite{KU}. It now suffices to take $a:=\lceil\log_pn^2\rceil$ and to note that $\log p$ is dominated by $2\log n$. Let $E(Y)$ be a monic lift of $\bar E(Y)$, then the ring $\Z_p[Y]/E(Y)$ has at least $n^2$ units and we conclude the proof with Theorem \ref{thmKU}.\hfill$\square$

We note that the use of the exponent $1+\epsilon$ in all our complexity estimates has the classical meaning that for every $\epsilon>0$ an algorithm exists with this estimate. For most results only logarithmic factors are needed (e.g.\ $\O\left(n\log n\right)$ instead of $\O\left(n^{1+\epsilon}\right)$), but we choose a more uniform formulation.

\subsection{Elementary operations: $+$, $-$, $\cdot$, $/$}\label{ssBasic}
Proposition \ref{propZpAdd} is essentially Corollary 11.10 in \cite{MCA}.
\begin{proposition}\label{propZpAdd}
Let $\alpha,\beta\in\Z_p\bmod p^N$. We can compute $\alpha+\beta$ and $\alpha-\beta$ in time $\O\left(N\log p\right)$. We can compute $\alpha\cdot \beta$ and if $\beta$ is a unit also $1/\beta$ in time $\O\left( (N\log p)^{1+\epsilon}\right)$.\end{proposition}
\begin{proposition}\label{propZqAdd}
Let $\alpha,\beta\in\Z_{p^n}\bmod p^N$. We can compute $\alpha+\beta$ and $\alpha-\beta$ in time $\O\left(Nn\log p\right)$. We can compute $\alpha\cdot \beta$ and if $\beta$ is a unit also $1/\beta$ in time $\O\left( (Nn\log p)^{1+\epsilon}\right)$.\end{proposition}
\textsc{Proof.} The result for $\alpha\pm\beta$ is entirely straightforward. Corollary 9.7 in \cite{MCA} implies that $\alpha\cdot\beta$ can be computed in the same amount of time --- up to a constant --- as the product of two polynomials of degree $n$ in $\Z_p[X]\bmod p^N$, which is $\O\left( (Nn\log p)^{1+\epsilon}\right)$ by Theorem 8.23 in the same book and Proposition \ref{propZpAdd} above. Over $\F_{p^n}$ the inverse of the reduction $\bar\beta$ of $\beta$ modulo $p$ can be computed in time $\O\left( (n\log p)^{1+\epsilon}\right)$ by Corollary 11.6 in \cite{MCA}. Afterwards $1/\bar\beta$ can be Newton lifted to $1/\beta$ by Algorithm 9.10 in \cite{MCA} which clearly has the required asymptotic complexity.\hfill$\square$

\subsection{Root finding (Newton iteration)}\label{ssNewton}
In this section we assume that some root, which is not a multiple root, is already known modulo $p$.
\begin{proposition}\label{propNewtonZp}
Let $f(Y)$ be a polynomial over $\Z_p\bmod p^N$ of degree $m$, and $y_0\in\Z_{p^n}\bmod p^N$ such that $f(y_0)\equiv 0\bmod p$ and $\frac{\partial f}{\partial Y}(y_0)\not\equiv 0\bmod p$. Then we can compute $y\in\Z_{p^n}\bmod p^N$ such that $y\equiv y_0\bmod p$ and $f(y)\equiv 0\bmod p^N$ in time $\O\left( (N(n+m)\log p)^{1+\epsilon}\right)$.\end{proposition}
\begin{proposition}\label{propNewtonZq}
Let $f(Y)$ be a polynomial over $\Z_{p^n}\bmod p^N$ of degree $m$, and $y_0\in\Z_{p^n}\bmod p^N$ such that $f(y_0)\equiv 0\bmod p$ and $\frac{\partial f}{\partial Y}(y_0)\not\equiv 0\bmod p$. Then we can compute $y\in\Z_{p^n}\bmod p^N$ such that $y\equiv y_0\bmod p$ and $f(y)\equiv 0\bmod p^N$ in time $\O\left( (Nnm\log p)^{1+\epsilon}\right)$.\end{proposition}
\textsc{Proof.} Both propositions can easily be proven by using classical $p$-adic Newton iteration with quadratic convergence, Algorithm 9.22 in \cite{MCA}. Note that for Proposition \ref{propNewtonZq} we need the obvious generalization of Theorem \ref{thmKU1} to polynomials over $\Z_{p^n}\bmod p^N$. The complexity estimates are entirely straightforward.\hfill$\square$

\subsection{Galois conjugates}\label{ssGalois}
We denote with $\sigma$ the $p$-th power Frobenius automorphism on $\Z_{p^n}$. Note that $\sigma^n$ is the identity map. \begin{proposition}\label{propSigma}
Let $\alpha\in\Z_{p^n}\bmod p^N$ and $0<k<n$ an integer. We can compute $\sigma^k(\alpha)$ in time $\O\left( ((N+\log p)n\log p)^{1+\epsilon}\right)$.
\end{proposition}
\textsc{Proof.} Let $\F_{p^n}\cong \F_p[\bar x]/\bar\varphi(\bar x)$ be the `reduction modulo $p$' of $\Z_{p^n}$, and $\bar\sigma$ the $p$-th power Frobenius on it. Clearly we can compute $\bar\sigma(\bar x)=\bar x^p$ in $\F_{p^n}$ in time $\O\left( (n\log^2 p)^{1+\epsilon}\right)$. In order to compute $\bar\sigma^k(\bar x)=\bar x^{p^k}$ we use the following lemma.
\begin{lemma}\label{lemCompFq} Given the polynomials $A(\bar x) := (\bar x^{p^a}\bmod \bar\varphi(\bar x))$ and $B(\bar x) := (\bar x^{p^b}\bmod \bar\varphi(\bar x))$ for some integers $a,b\geq 1$, we have that $A(B(\bar x))\equiv\bar x^{p^{a+b}}\bmod \bar\varphi(\bar x)$, and this composition can be computed in time $\O\left( (n\log p)^{1+\epsilon}\right)$.\end{lemma}
\textsc{Proof of the lemma.} It is easy to verify that $A(B(\bar x))\bmod \bar\varphi(\bar x) = \bar x^{p^{a+b}}\bmod \bar\varphi(\bar x)$, using the fact that $B(\bar x)$ is a root of $\bar\varphi(\bar x)$. Now Theorem \ref{thmKU1} (for $N=1$) gives the lemma.\hfill$\square$

\noindent
\textsc{Proof of the proposition (continued).} The idea to compute $\bar\sigma^k(\bar x)$ is to use the binary representation of $k$ combined with the lemma. The general algorithm is similar to the classical repeated squaring technique (Algorithm 4.8 in \cite{MCA}), we explain here only the easier case where $k=2^m$ for an integer $m\geq 1$. 
The procedure is quite obvious: compute recursively $A_i(\bar x)=A_{i-1}(A_{i-1}(\bar x))\bmod \varphi(\bar(x))$ with $A_0(\bar x)=\sigma(\bar x)$. Lemma \ref{lemCompFq} yields $A_1(\bar x)=\bar x^{p^2}\bmod \bar\varphi(\bar x)$, $A_2(\bar x)=\bar x^{p^4}\bmod \bar\varphi(\bar x)$, \ldots, $A_m(\bar x)=\bar x^{p^{2^m}}\bmod \bar\varphi(\bar x)$. Only $m=\log_2 k\leq \log n$ steps are required, hence if we know $\bar \sigma(\bar x)$, we can compute $\bar\sigma^k(\bar x)$ in time $\O\left( \log k(n\log p)^{1+\epsilon}\right)=\O\left( (n\log p)^{1+\epsilon}\right)$.

Because $\sigma^k(x)$ is a root of $\varphi(X)$ and $\bar\varphi(X)$ is squarefree, we can now apply Proposition \ref{propNewtonZp} in order to lift $\bar\sigma^k(\bar x)$ to $\sigma^k(x)$ modulo $p^N$ in time $\O\left( (Nn\log p)^{1+\epsilon}\right)$. For $\alpha(x)\in\Z_p[x]/\varphi(x)$ we have $\sigma^k(\alpha(x))=\alpha(\sigma^k(x))\bmod \varphi(x)$, and hence Theorem \ref{thmKU1} allows us to compute this last expression with precision $p^N$ in time $\O\left( (Nn\log p)^{1+\epsilon}\right)$, thereby proving the proposition.\hfill$\square$

\begin{corollary}\label{corFq}
Let $\alpha\in\F_{p^n}$, $\bar\sigma$ the Frobenius automorphism and $0<k<n$. Then we can compute $\bar\sigma^k(\alpha)$ in time $\O\left((n\log^2p)^{1+\epsilon}\right)$.\end{corollary}
\textsc{Proof.} With $\F_{p^n}$ given as $\F_p[\bar x]/\bar\varphi(\bar x)$, we have shown above that $\bar\sigma^k(\bar x)$ can be computed in time $\O\left((n\log^2p)^{1+\epsilon}\right)$. Now writing $\alpha$ as $\alpha(\bar x)$ gives $\bar\sigma^k(\alpha(\bar x))=\alpha(\bar\sigma^k(\bar x))$, hence Theorem \ref{thmKU1} gives the Corollary.\hfill$\square$

\subsection{Equations with Frobenius}\label{ssEquations}
In this section we merely rephrase results from \cite{HEHCC} using faster Frobenius computations.
\begin{proposition}\label{propASeq}
Let $\alpha,\beta,\gamma\in\Z_{p^n}\bmod p^N$ with $\beta\equiv 0\bmod p$. We can compute the (unique)  solution $X$ in $\Z_{p^n}\bmod p^N$ of $\alpha \sigma(X)+\beta X+\gamma=0$ in time $\O\left( ((N+\log p)n\log p)^{1+\epsilon}\right)$.\end{proposition}
\begin{proposition}\label{propGenNewton}
Let $\phi(Y,Z)$ be a polynomial over $\Z_{p^n}\bmod p^N$ for which the evaluation of $\phi$, $\partial\phi/\partial Y$ and $\partial\phi/\partial Z$ in any $(\alpha, \beta)\in(\Z_{p^n}\bmod p^N)^2$ requires at most $\psi$ arithmetic operations in $\Z_{p^n}\bmod p^N$. Suppose we have $x_0\in\Z_{p^n}\bmod p^N$ such that $\phi(x_0,\sigma(x_0))\equiv 0\bmod p^{2k+1}$ with $k:=\text{ord}_p(\frac{\partial \phi}{\partial Z}(x_0,\sigma(x_0))$. Then we can compute $X\in\Z_{p^n}\bmod p^{N+k}$ such that $\phi(X,\sigma(X))\equiv 0\bmod p^{N+k}$ and $X\equiv x_0\bmod p^{k+1}$ in time $\O\left( ((\psi N+\log p)n\log p)^{1+\epsilon}\right)$.\end{proposition}
\textsc{Proof of Proposition \ref{propASeq}.} In Section 12.6.1 of \cite{HEHCC}, an algorithm by Lercier and Lubicz \cite{LeLu} is explained that computes the solution $X$. Its complexity is determined by Algorithm 12.18 of \cite{HEHCC}, which gives $\O\left(\log n((N+\log p)n\log p)^{1+\epsilon}\right)$ if we use Proposition \ref{propSigma} above.\hfill$\square$

\noindent\textsc{Proof of Proposition \ref{propGenNewton}.} Again we recycle an algorithm of \cite{HEHCC}, namely Algorithm 12.23 which computes a generalized Newton lift. Except for $\O\left(\log N\right)$ times an evaluation of $\phi$, $\partial\phi/\partial Y$ and $\partial\phi/\partial Z$, its complexity it the same as the one given in Proposition \ref{propASeq} above. Hence the total complexity is bounded by
$\O\left( \psi\log N(Nn\log p)^{1+\epsilon} + ((N+\log p)n\log p)^{1+\epsilon}\right)$.\hfill$\square$

\subsection{Teichm\"{u}ller lift}\label{ssTeichLift}
\begin{proposition}\label{propTeichLift}
Given $\alpha\in\Z_{p^n}\bmod p$, we can compute the Teichm\"{u}ller lift of $(\alpha\bmod p)$ in time $\O\left( (Nn\log^2 p)^{1+\epsilon}\right)$.
\end{proposition}
\textsc{Proof.} As pointed out in Section 12.8.1 of \cite{HEHCC}, we can use Proposition \ref{propGenNewton} for the polynomial $\phi(Y,Z)=Y^p-Z$ with $x_0=\alpha$ and $k=0$. Evaluating $\phi$, $\partial \phi/\partial Y$ and $\partial \phi/\partial Z$ requires $\O\left(\log p\right)$ elementary operations in $\Z_{p^n}\bmod p^N$ and we find the proposition.\hfill$\square$

\subsection{Minimal polynomial, trace and norm}\label{ssMinPol}
\begin{proposition}\label{propMinPol}
Let $\alpha\in\Z_{p^n}\bmod p^N$. We can compute the minimal polynomial modulo $p^N$ of $\alpha$ over $\Z_p$ in time $\O\left( (Nn\log p)^{1+\epsilon}\right)$.
\end{proposition}
\begin{corollary}\label{corTraceNorm}
Let $\alpha\in\Z_{p^n}\bmod p^N$. We can compute the trace $\textnormal{Tr}(\alpha)$ and norm $\textnormal{N}(\alpha)$ over $\Z_p\bmod p^N$ in time $\O\left( (Nn\log p)^{1+\epsilon}\right)$.\end{corollary}
A \emph{Teichm\"{u}ller modulus} is the minimal polynomial of a Teichm\"{u}ller lift (see Section 12.1 of \cite{HEHCC}), or equivalently a divisor of $X^{p^n}-X$ for appropriate $n$.
\begin{corollary}\label{corTeichMod}
Given $\F_{p^n}\cong\F_p[x]/\bar\varphi(x)$, we can compute a Teichm\"{u}ller modulus $F(X)$ modulo $p^N$ which equals $\bar\varphi(X)$ modulo $p$ in time $\O\left( (Nn\log^2 p)^{1+\epsilon}\right)$.
\end{corollary}
\textsc{Proof of proposition \ref{propMinPol}.}
We follow an idea of \cite{Rifa} and \cite{Ly} as explained in Section 3 of \cite{ShFastConstr}. Define the linear operator $P:\Z_p[x]/\varphi(x)\to\Z_p$ by $P(1) := 1$ and $P(x)=P(x^2)=\ldots=P(x^{n-1})=0$. We can compute --- using the fast modular power projection of Theorem 7.7 in \cite{KU} --- the sequence $P(1)$, $P(\alpha), \ldots, P(\alpha^{2n-1})$ in essentially linear time $\O\left((Nn\log p)^{1+\epsilon}\right)$. The minimal polynomial $c(X)$ of $\{P(\alpha^i)\}_{i\geq 0}$ equals the minimal polynomial of $\alpha$ modulo $p^N$, and Step 2 of Shoup's algorithm refers to the fact that one can obtain this minimal polynomial from the (fast) extended Euclidean algorithm for $$g(X)=\sum_{i=0}^{2n-1}P(\alpha^i)X^{2n-1-i}\qquad\textnormal{and}\qquad f(X)=X^{2n}.$$
Indeed, knowing a Euclidean expansion $c(X)g(X)+q(X)f(X)=r(X)$ for some remainder $r(X)$ of degree at most $n-1$ and with $c(X)$ of minimal degree, implies that $c(X)$ is the minimal polynomial of $\alpha$.
\hfill$\square$

We note that for computing $\textnormal{N}(\alpha)$ a much more elegant algorithm was given by Harley, see Section 12.8.5.c in \cite{HEHCC}. Namely, if we write $\alpha$ as $\alpha(x)$, the resultant formula $\textrm{N}(\alpha)=\textrm{Res}_X(\varphi(X),\alpha(X))$ can be computed in the same amount of time as in Corollary \ref{corTraceNorm}, using a variant of Moenck's extended gcd algorithm \cite{Moenck}.

\end{document}